\documentclass[11pt]{article}

\usepackage{amssymb}
\usepackage{latexsym}
\usepackage{amsmath}[1996/11/01]

\newcommand{\eeq}{\end{equation}}
\newcommand{\beql}[1]{\begin{equation}\label{#1}}
\newcommand{\beq}{\begin{equation}}
\newcommand{\eqn}[1]{(\ref{#1})}
\newcommand{\ZZ}{{\mathbb Z}}

\newcommand{\QQ}{{\mathbb Q}}
\newcommand{\FF}{{\mathbb F}}
\newcommand{\CC}{{\mathbb C}}
\newcommand{\sC}{{\mathcal C}}
\DeclareMathOperator{\Mat}{Mat}
\DeclareMathOperator{\Tr}{Tr}
\DeclareMathOperator{\Sym}{Sym}

\newtheorem{theorem}{Theorem}

\makeatletter
\def\@sect#1#2#3#4#5#6[#7]#8{\ifnum #2>\c@secnumdepth
     \def\@svsec{}\else
     \refstepcounter{#1}\edef\@svsec{\csname the#1\endcsname.\hskip .75em }\fi
     \@tempskipa #5\relax
      \ifdim \@tempskipa>\z@
        \begingroup #6\relax
          \@hangfrom{\hskip #3\relax\@svsec}{\interlinepenalty \@M #8\par}%
        \endgroup
       \csname #1mark\endcsname{#7}\addcontentsline
         {toc}{#1}{\ifnum #2>\c@secnumdepth \else
                      \protect\numberline{\csname the#1\endcsname}\fi
                    #7}\else
        \def\@svsechd{#6\hskip #3\@svsec #8\csname #1mark\endcsname
                      {#7}\addcontentsline
                           {toc}{#1}{\ifnum #2>\c@secnumdepth \else
                             \protect\numberline{\csname the#1\endcsname}\fi
                       #7}}\fi
     \@xsect{#5}}
\def\@begintheorem#1#2{\it \trivlist \item[\hskip \labelsep{\bf #1\ #2.}]}
\makeatother

\begin{document}

\begin{center}
{\Large\bf Gleason's Theorem on Self-Dual Codes and Its Generalizations} \\

\vspace*{+.2in}
N. J. A. Sloane \\ [+.2in]
AT\&T Shannon Labs, Florham Park, NJ, USA \\
\vspace*{+.2in}
September 30, 2006; corrected October 17, 2006 \\
\vspace*{+.2in}
TO EIICHI BANNAI, ON THE OCCASION OF HIS $60$TH BIRTHDAY. \\
\vspace*{+.2in}
{\bf Abstract}
\vspace*{+.2in}
\end{center}

One of the most remarkable theorems in coding theory is Gleason's 1970
theorem about the weight enumerators of self-dual codes.
In the past
36 years there have been hundreds of papers written about generalizations
and
applications of this theorem to different types of codes, always
on a case-by-case basis.
In this talk I will state the theorem and then
describe the far-reaching generalization that Gabriele Nebe, Eric Rains and I have
developed which includes all the earlier generalizations at once.
The full proof has just appeared in our book {\em Self-Dual Codes
and Invariant Theory} (Springer, 2006).

This paper is based on my talk at
the conference on Algebraic Combinatorics in honor of
Eiichi Bannai, held in Sendai, Japan, June 26--30, 2006.

\setlength{\baselineskip}{1.5\baselineskip}
\section{Motivation}

Self-dual codes are important
because they intersect with
\begin{itemize}
\item
communications
\item
combinatorics
\item
block designs, spherical designs
\item
group theory
\item
number theory
\item
sphere packing
\item
quantum codes
\item
conformal field theory, string theory
\end{itemize}

\section{Introduction}

In classical coding theory (as for example in MacWilliams and Sloane
\cite{MS77}), a {\em code} $C$ of length $N$ over a field $\FF$ is a subspace of $\FF^N$.
The {\em dual code} is
$$C^\perp := \{ u \in \FF^N : ~u \cdot c =0 , ~~\forall c \in C \} \,.$$

\paragraph{Example:}
 $C = \{ 000, 111\}$, $C^\perp = \{000, 011, 101, 110\}$ with $\FF = \FF_2$.
The weight enumerators of these two codes are
$$W_C (x,y) = x^3 + y^3 , ~~
W_{C^\perp} (x,y) = x^3 + 3xy^2 \,.
$$
A code is {\em self-dual} if $C= C^\perp$.
For example, the binary code $i_2 := \{00, 11\}$ is self-dual, with weight enumerator
\beql{Eq1}
W_{i_2} (x,y) = x^2 + y^2 \,.
\eeq

\paragraph{Example:}
The Hamming code $h_8$ of length 8 is self-dual.
This is the binary code with generator matrix:
$$
\begin{array}{|cccccccc|}
\multicolumn{1}{c}{\infty} & 0 & 1 & 2 & 3 & 4 & 5 & \multicolumn{1}{c}{6} \\ \hline
1 & 1 & 1 & 1 & 1 & 1 & 1 & 1 \\
0 & 1 & 1 & 1 & 0 & 1 & 0 & 0 \\
0 & 0 & 1 & 1 & 1 & 0 & 1 & 0 \\
0 & 0 & 0 & 1 & 1 & 1 & 0 & 1 \\ \hline
\end{array}
$$
The second row of the matrix has 1's under the quadratic residues 0, 1, 2 and 4 $\bmod~7$.
The remaining rows are obtained by fixing the infinity coordinate and cycling the other coordinates.
This code has weight enumerator
\beql{Eq2}
W_{h_8} (x,y) = x^8 + 14 x^4 y^4 + y^8 \,.
\eeq
As can be seen from the generator matrix, this code is closely related to the 
incidence matrix of the projective plane of order 2.

If we replace the prime 7 in this construction by 23, we get the binary 
Golay self-dual code $g_{24}$ of length 24, with generator matrix as follows:
$$
\begin{array}{|cccccccccccccccccccccccc|} \hline
1, & 1, & 1, & 1, & 1, & 1, & 1, & 1, & 1, & 1, & 1, & 1, & 1, & 1, & 1, & 1, & 1, & 1, & 1, & 1, & 1, & 1, & 1, & 1\\
0, & 1, & 1, & 1, & 1, & 1, & 0, & 1, & 0, & 1, & 1, & 0, & 0, & 1, & 1, & 0, & 0, & 1, & 0, & 1, & 0, & 0, & 0, & 0\\
0, & 0, & 1, & 1, & 1, & 1, & 1, & 0, & 1, & 0, & 1, & 1, & 0, & 0, & 1, & 1, & 0, & 0, & 1, & 0, & 1, & 0, & 0, & 0\\
0, & 0, & 0, & 1, & 1, & 1, & 1, & 1, & 0, & 1, & 0, & 1, & 1, & 0, & 0, & 1, & 1, & 0, & 0, & 1, & 0, & 1, & 0, & 0\\
0, & 0, & 0, & 0, & 1, & 1, & 1, & 1, & 1, & 0, & 1, & 0, & 1, & 1, & 0, & 0, & 1, & 1, & 0, & 0, & 1, & 0, & 1, & 0\\
0, & 0, & 0, & 0, & 0, & 1, & 1, & 1, & 1, & 1, & 0, & 1, & 0, & 1, & 1, & 0, & 0, & 1, & 1, & 0, & 0, & 1, & 0, & 1\\
0, & 1, & 0, & 0, & 0, & 0, & 1, & 1, & 1, & 1, & 1, & 0, & 1, & 0, & 1, & 1, & 0, & 0, & 1, & 1, & 0, & 0, & 1, & 0\\
0, & 0, & 1, & 0, & 0, & 0, & 0, & 1, & 1, & 1, & 1, & 1, & 0, & 1, & 0, & 1, & 1, & 0, & 0, & 1, & 1, & 0, & 0, & 1\\
0, & 1, & 0, & 1, & 0, & 0, & 0, & 0, & 1, & 1, & 1, & 1, & 1, & 0, & 1, & 0, & 1, & 1, & 0, & 0, & 1, & 1, & 0, & 0\\
0, & 0, & 1, & 0, & 1, & 0, & 0, & 0, & 0, & 1, & 1, & 1, & 1, & 1, & 0, & 1, & 0, & 1, & 1, & 0, & 0, & 1, & 1, & 0\\
0, & 0, & 0, & 1, & 0, & 1, & 0, & 0, & 0, & 0, & 1, & 1, & 1, & 1, & 1, & 0, & 1, & 0, & 1, & 1, & 0, & 0, & 1, & 1\\
0, & 1, & 0, & 0, & 1, & 0, & 1, & 0, & 0, & 0, & 0, & 1, & 1, & 1, & 1, & 1, & 0, & 1, & 0, & 1, & 1, & 0, & 0, & 1\\ \hline
\end{array}
$$

\noindent
This has weight enumerator
\beql{Eq3}
W_{g_{24}}(x,y) ~=~ x^{24} + 759 x^{16} y^8 + 2576 x^{12} y^{12} + 759 x^8 y^{16} + y^{24} \,.
\eeq

\section{MacWilliams' Theorem, 1962}

In her Ph.D. thesis at Harvard in 1962, Jessie MacWilliams \cite{MacW63}
 showed that the weight enumerator of the dual of
a linear code is determined by the weight enumerator
of the code:

\begin{theorem}\label{th1}
For a code $C$ over $\FF_q$,
\beql{Eq4}
W_{C^\perp} (x,y) = \frac{1}{|C|} W_C (x+ (q-1) y , x-y ) \,.
\eeq
\end{theorem}
The proof uses the Poisson summation formula,
 in the form that says that the sum of a function $f$ over a vector space is equal to the average of the appropriate Fourier transform of $f$ over the dual vector space.

\paragraph{Corollary.}
If $C$ is self-dual, $W_C (x,y)$ is fixed under
the ``MacWilliams'' transformation
\beql{Eq5}
\left( \begin{array}{c}
x \\ y
\end{array}\right) 
\mapsto 
\frac{1}{\sqrt{q}}
\left(
\begin{array}{cc}
1 & q-1 \\ 1 & -1
\end{array}
\right)
\left( \begin{array}{c}
x \\ y
\end{array}\right) 
\,.
\eeq

\section{First there were four types}

As can be seen from \eqn{Eq2} and \eqn{Eq3},
for some binary self-dual codes the Hamming weights
of all the codewords (the powers of $y$) are multiples of 4.
In other cases (as in \eqn{Eq1}) the weights may only be even.
Gleason and Pierce showed that there are essentially only four possibilities for this phenomenon to occur with self-dual codes over fields:
\begin{theorem}\label{th2}
{\rm (Gleason-Pierce (1967), see \cite{AMT67}, \cite{Me64}).}
If $C$ is a self-dual code over $\FF_q$ with Hamming weights divisible by $m$,
 then one of the following holds:
\begin{itemize}
\item[{\rm I})]
$q=2$ $(\Rightarrow m =2 )$
\item[{\rm II})]
$q=2$ and $m=4$
\item[{\rm III})]
$q=3$ $(\Rightarrow m=3)$
\item[{\rm IV})]
$q=4$ and Hermitian $(\Rightarrow m=2)$ or $q=4$ and Euclidean,
\item[{}]
or else $c=2$, $q$ is
arbitrary, $N$ is even and $W(x,y) = (x^2 + (q-1) y^2)^{N/2}$.
\end{itemize}
\end{theorem}

Because of this theorem, self-dual codes falling into one of those four classes came to be known as codes of Types I, II, III and IV, respectively.

Incidentally, the codes with $c=2$ mentioned in the final clause of the theorem still have not been fully classified (see \cite{Me64}).

\section{Gleason's Theorem (1970, Nice)}

At the International Congress of Mathematicians in Nice, 1970, Gleason established the following result.

\begin{theorem}\label{th3}
{\rm (Gleason \cite{Gle70}).}
If $C$ is a self-dual
code of one of the four types mentioned in Theorem \ref{th2} then the weight enumerator of $C$ belongs to the polynomial ring $\CC [f, g]$, where:
\end{theorem}
\renewcommand{\arraystretch}{1.4}
\begin{center}
\begin{tabular}{|c|c|c|}
\hline
Type & $f$ & $g$ \\ \hline
I & $x^2 + y^2$ & $x^2 y^2 (x^2 - y^2 )^2$ \\
~ & $i_2$ & Hamming code $h_8$ \\ \hline
II & $x^8 + 14 x^4 y^4 + y^8$ & $x^4 y^4 (x^4 - y^4)^4$ \\
~ & Hamming code $h_8$ & binary Golay code $g_{24}$ \\ \hline
III & $x^4 + 8xy^3$ & $y^3 (x^3 - y^3)^3$ \\
& tetracode & ternary Golay code \\ \hline
IV* & $x^2 + 3y^2$ & $y^2 (x^2 - y^2 )^2$ \\
& $i_2 \otimes \FF_4$ & hexacode \\ \hline
\end{tabular}
\end{center}
{\footnotesize {*In fact Gleason omitted this case, which was first given in {\rm \cite{MMS72}}.}}
\renewcommand{\arraystretch}{1.0}

Under each polynomial we have written the name of 
a code whose weight enumerator leads to that polynomial.
For example, the theorem states that the weight enumerator of a Type I self-dual code belongs to
the ring generated by the weight enumerators of the codes $i_2$ and $h_8$,
that is, by $f=W_{i_2} =x^2 + y^2$ (see \eqn{Eq1}) and $W_{h_8}$ (see \eqn{Eq2}).
It is simpler to replace $W_{h_8}$ as a generator of this ring by
$$g := \frac{1}{4} \left( f^4 - W_{h_8} \right) =
x^2 y^2 (x^2 - y^2 ) \,,$$
as in the first row of the table.

In the following years many generalizations of Gleason's theorem were published, for example to self-dual codes over other fields $(\FF_5 , \ldots )$, to biweight enumerators, split weight enumerators, codes containing the all-ones
vector, etc.

The main applications of these theorems are
in the classification of 
self-dual codes of moderate lengths,
and
in the determination of the
optimal (or {\em extremal}) codes of the various Types.
The book \cite{cliff2}
contains an extensive bibliography.

\section{$\ZZ / 4 \ZZ$ appears!}

In the early 1990's, coding theory changed forever when it was discovered that
certain infamous nonlinear binary codes were really {\em linear} (and in some cases self-dual) codes over the ring $\ZZ / 4 \ZZ$ of integers $\bmod~4$.
For example, the Nordstrom-Robinson code is a famous nonlinear binary code of
length 16 that contains 256 codewords and has minimal distance 6, more than is possible with any linear code of the same length and the same number of
codewords (cf. \cite{MS77}).
Although nonlinear, its weight enumerator behaves like that of a linear
binary code --- it is fixed under \eqn{Eq5} (with $q=2$).
In 1992, Forney, Trott and I \cite{Me183} showed that this code is really a linear
self-dual code over $\ZZ / 4 \ZZ$, a code already known as the {\em octacode}
(cf. \cite{Me168}, \cite{SPLAG}).

This work was extended by Hammons, Kumar, Calderbank, Sol\'{e} and me in \cite{Me184} to cover
many other families of binary nonlinear codes, and this in turn was followed
by numerous other papers that studied self-dual codes over the rings $\ZZ / m \ZZ$ for integers $m \ge 4$ (again see \cite{cliff2} for references).

\section{And then there were nine!}

In 1998, Eric Rains and I wrote a 120-page survey for the {\em Handbook of Coding Theory} \cite{RS98} in which we distinguished nine types of self-dual codes,
extending the original four types to include such families
as linear codes over $\ZZ / 4 \ZZ$, linear codes over
$\ZZ / m \ZZ$ for $m \ge 4$, additive codes over $\FF_4$, etc.
Again each version of Gleason's theorem was treated separately.

\section{Higher-genus weight enumerators}

In the mid-1990's there was a major breakthrough.
As the result of a really amazing coincidence, we were led to investigate a certain family of
``Clifford groups'' $\sC (m)$ with structure $2_+^{1+2m} \cdot O_{2m}^+ (2)$.
The story of this astonishing coincidence can be found in \cite{Me223}
and \cite{cliff2}, so I will not repeat it here.
Studying these ``Clifford''  groups led to breakthroughs in quantum codes
\cite{Me223} and to generalizations of Gleason's theorem to higher-genus or multiple weight enumerators \cite{cliff1}.

\section{The new book}

After writing \cite{cliff1}, we realized that the arguments used to handle the invariants of the Clifford groups could be extended to handle other classes of
self-dual codes.
The result is a far-reaching generalization of Gleason's theorem which defines the ``Type'' of a self-dual code in such a way that the weight enumerator of any code of that Type belongs to the invariant ring of a certain ``Clifford-Weil'' group associated with the Type, and furthermore that this invariant ring is spanned by weight enumerators of codes of that Type.

These are the two properties that previously had to be proved for each Type on a case-by-case basis.
Now we know that this is automatically true, provided the codes fall into
one of certain very general classes.

The proof of the general theorem is not easy, and occupies perhaps 150 pages of the new 400-page book, ``Self-Dual Codes and Invariant Theory'' \cite{cliff2}.

For me, the book represents the culmination of thirty-five years of work.

In the rest of this talk I will give an outline of our approach, omitting all the technical details (and the category theory).

\section{Notation: codes over rings}

We will use the following notation:
$$
\begin{array}{lll}
R & = & \mbox{ground ring $=$ ring with unit} \\ [+.1in]
V & = & \mbox{left $R$-module $=$ alphabet} \\ 
&& \mbox{(usually we assume $R$ and $V$ are finite)} \\ [+.1in]
C & = & \mbox{code of length $N$} \\ [+.1in]
& = & \mbox{$R$-submodule of $V^N$} \\ [+.1in]
\multicolumn{3}{l}{c \in C, ~ r \in R \Rightarrow rc \in C\,.}
\end{array}
$$

\paragraph{Dual code:}
$$
\begin{array}{c}
\beta = \mbox{nonsingular bilinear form}\,, ~~
\beta : ~ V \times V \to \QQ / \ZZ \\ [+.15in]
\mbox{E.g.} ~~\beta (u,v) = \frac{1}{2} uv~({\rm in~the~binary~case}) \\ [+.15in]
C^\perp := \left\{ u \in V^N : ~\displaystyle\sum_{i=1}^N \beta (u_i, c_i) =0, ~ \forall c \in C \right\}
\end{array}
$$

\section{Weight enumerators}

Let $C \le V^N$ be a code, where the alphabet $V= \{ v_0 , v_1, \ldots \}$.
The {\em complete weight enumerator} of $C$ is:
$$\mbox{cwe} (C) ~:=~
\sum_{c \in C} ~\prod_{i=1}^N x_{c_i}~\in~\CC [x_{v_0}, x_{v_1} , \ldots ] \,.$$

\paragraph{Example:}
$\FF_4 = \{0,1, \omega, \bar{\omega} \}$,
$$\mbox{codeword} ~~ c=0011\omega\bar{\omega}~ \Rightarrow ~ x_0^2 x_1^2 x_\omega x_{\bar{\omega}} \,.$$
The {\em symmetrized weight enumerator} is obtained by identifying $x_v$ and $x_w$ in $\mbox{cwe} (C)$ if we do not need to distinguish $v$ and $w$.
E.g. we usually set $x_{\bar{\omega}} = x_\omega$ for codes over $\FF_4$.
The {\em Hamming weight enumerator} is obtained from $\mbox{cwe} (C)$ by setting
$x_0 =x$ and all other $x_v = y$.

\section{Biweight or genus-2 weight enumerator}

Take an ordered pair of codewords $b$, $c \in C$ in all possible ways
and write one above the other:
$$
\left[ \begin{array}{c}
b \\ c
\end{array}\right]
 = \left[
\begin{array}{ccccccc}
0 & 0 & 1 & 1 & \omega & \bar{\omega} & \cdots \\
0 & 1 & 0 & 1 & \omega & \bar{\omega} & \cdots
\end{array}
\right] \,.
$$
Then the {\em biweight} or {\em genus}-2 weight enumerator
of $C$ is
$$\mbox{cwe}_2 (C) := \sum_{(b,c) \in C \times C}~ \prod_{i=1}^N
x_{\left( \begin{array}{c} b_i \\ c_i \end{array}\right)} \,.
$$

\paragraph{Remark:}
$$\mbox{cwe}_2 (C) = \mbox{cwe} (C \otimes R^2 ) \,.$$
For we have
$$C \otimes R^2 \le V^N \otimes R^2 \cong V^{2N} \cong (V^2)^N \,.$$
Note that the ground ring for $C \otimes R^2$ is
$$\Mat_2 (R) \,,$$
the ring of $2 \times 2$ matrices over $R$.
So, even in the case of classical binary codes, we need to
use noncommutative rings when we consider higher-genus weight enumerators!

\section{Extra conditions}

Often one wants to consider self-dual codes with certain additional properties,
for example that the weights are divisible by 4, or the code contains
the all-ones vector.
Some of these properties can be included in the new notion of Type, provided they can be described in terms of ``quadratic mappings''.
Oversimplifying (see \cite[Chapter 1]{cliff2} for the precise definition), a
quadratic mapping is a map from $V$ to $\QQ / \ZZ$ which is the sum of a quadratic part and a linear part.
If $\Phi$ is a collection of quadratic mappings then we say that a code $C$ is
{\em isotropic} with respect to $\Phi$ if
$$\sum_{i=1}^N \phi (c_i ) =0 , \qquad \forall ~c \in C , ~~\forall \phi \in \Phi \,.
$$
\paragraph{Examples:}
\begin{itemize}
\item
$\phi (x) = \frac{1}{4} x^2$ (to get weights divisible by 4 in the binary case)
\item
$\phi (x) = \frac{1}{p} x$, $p$ odd (to ensure that $1 \in C$)
\item
$\phi (x) = \beta (x,x)$ (specialization of $\beta$, always present)
\end{itemize}

\section{The new definition of Type}

We say that a code $C \le V^N$ has
$${\rm Type} ~ \rho := (R, V, \beta, \Phi )$$
if $C$ is self-dual with respect to $\beta$ and isotropic with respect to $\Phi$.

\vspace*{+.1in}
\noindent{\bf Memo:}
Many details have been concealed here.
See \cite{cliff2} for further information.

We call $(R, V, \beta, \Phi )$ a {\em form ring}, adapting a term from algebraic $K$-theory (cf. Bak \cite{Bak1}).

\section{Symmetric idempotents}

A {\em symmetric idempotent} $\iota \in R^\ast$ satisfies $\iota^2 = \iota$ together with certain extra
conditions (see \cite{cliff2}), and has the property that there are ``left'' and ``right'' elements $l_\iota$ and $r_\iota$ associated with it such that
$$\iota = l_\iota r_\iota$$

\paragraph{Examples:}
\begin{itemize}
\item
$R = \ZZ / 6 \ZZ$: $\iota = 3 = 3 \cdot 3$ or
$~\iota = 4 = 4 \cdot 4$

\item
$R = \Mat_m (R' ):$ $\iota = {\rm diag} \{1,0,0, \ldots, 0\}$
\end{itemize}

\section{The Clifford-Weil group $\sC (\rho )$}

We associate with the form ring $\rho = (R, V, \beta, \Phi )$ a certain finite subgroup of $GL_{|V|} ( \CC )$ that we call the Clifford-Weil group $\sC (\rho )$.
This generalizes the familiar group of order 192 generated by
$$\left[ \begin{array}{cc}
1 & 0 \\ 0 & i
\end{array} \right] , \qquad \frac{1}{\sqrt{2}} \left[ \begin{array}{cr}
1 & 1 \\ 1 & -1
\end{array}\right]
$$
that arises from Gleason's theorem for Type II (or doubly-even) binary codes,
and also generalizes the Clifford groups $\sC (m)$ mentioned above.
The generators for $\sC (\rho )$ are:
$$
\begin{array}{ll}
\rho (r) : & x_v \mapsto x_{rv} , \qquad \forall r \in R^\ast \quad (\mbox{because $C$ is a code}) \\ [+.1in]
\rho (\phi ): & x_v \mapsto e^{2 \pi i \phi (v)} x_v , \qquad
\forall \phi \in \Phi 
\quad \mbox{(because~$C$ is isotropic)}
\end{array}
$$
and the ``MacWilliams'' transformations (generalizing \eqn{Eq5}):
for every symmetric idempotent $\iota = l_\iota r_\iota$ the associated MacWilliams transformation is
$$h_{\iota , r} :
x_v \mapsto \frac{1}{\sqrt{| \iota V|}} \sum_{v \in V} e^{2 \pi i \beta (w, r_\iota v)} x_{w+ (1-\iota) v}\,.
$$

\paragraph{Example of $h_{\iota, r}$:}
$$
\begin{array}{l}
R = V = \ZZ / 6 \ZZ \,, \\ [+.1in]
\mbox{two symmetric idempotents $3=3 \cdot 3$, $4 = 4 \cdot 4$}
\end{array}
$$
For $\iota =3$, $r_{\iota} =3$, $|3V|=2$:
$$h_{3, r_3} = \frac{1}{\sqrt{2}} 
\left[ \begin{array}{cccccc}
+ & 0 & 0 & + & 0 & 0 \\
0 & - & 0 & 0 & + & 0 \\
0 & 0 & + & 0 & 0 & + \\
+ & 0 & 0 & - & 0 & 0 \\
0 & + & 0 & 0 & + & 0 \\
0 & 0 & + & 0 & 0 & -
\end{array}
\right] \,,
$$
where $+$ stands for $+1$ and $-$ for $-1$.

Our reasons for calling $\sC (\rho )$ the Clifford-Weil group are that (i)~when the groups $\sC (m)$ mentioned in \S8 act on the Barnes-Wall lattices,
they act as the full orthogonal group $O_{2m}^+ (2)$ on the Clifford algebra of the quadratic form, and
(ii)~in some situations $\sC (\rho )$ coincides with the groups
studied by Weil in his famous paper
``Sur certaines groups d'op\'{e}rateurs unitaires'' \cite{Weil64}.

\section{Quasi-chain rings}

Our main theorems will cover self-dual codes over a very large class of rings.

A {\em chain ring} is one in which the left ideals are linearly ordered by inclusion.

A {\em quasi-chain ring} is a direct product of matrix rings over chain rings.

\paragraph{Examples of quasi-chain rings:}
\begin{itemize}
\item
matrix rings over finite fields
\item
matrix rings over $\ZZ / m \ZZ$
\item
matrix rings over Galois rings
\end{itemize}
But not all rings are covered by the present theory.
Examples of rings that are not (yet) covered:
\begin{itemize}
\item
the group ring $\FF_3\Sym(3)$
\item
the matrix ring
$$
\left[ \begin{array}{cc}
\ZZ / 4 \ZZ & \ZZ / 4 \ZZ \\
2 \ZZ / 4 \ZZ & \ZZ / 4 \ZZ
\end{array}\right]
$$
\end{itemize}

\section{The main theorems}

\paragraph{Theorem.}
Let $R$ be a finite chain ring or quasi-chain ring, and let $\rho$ be the form ring
$$\rho = (R,V, \beta , \Phi ) \,.$$
Consider codes $C \le V^N$ of Type $\rho$.
Then (i) $\mbox{cwe} (C)$ belongs to the invariant ring $\mbox{Inv} ( \sC (\rho ))$, and (ii)~$\mbox{Inv} (\sC (\rho ))$ is spanned by the $\mbox{cwe} (C)$, where $C$ 
runs through codes of Type $\rho$.

The proof, as already mentioned, uses category theory and is long and hard, and takes up a good part of the book \cite{cliff2}.

We believe, but have not been able to prove, that the theorem should hold
without the restriction to quasi-chain rings.
We state this as the:

\paragraph{Weight Enumerator Conjecture:}
The theorem should hold for any finite ring $R$.

\section{An application}

In their 1999 paper ``Type II codes, even unimodular lattices and invariant rings''
\cite{BDHO99}, Bannai, Dougherty, Harada and Oura consider codes of
(in our new notation) Type $4_{{\rm II}}^{\ZZ}$.
The corresponding Clifford-Weil group has order 1536, and the ring to which the complete weight enumerators belong has Molien series
\beql{EqB1}
\frac{1+t^8 +2t^{16} + 2t^{24} + t^{32} + t^{40}}{(1-t^8 )^3 (1-t^{24})} \,.
\eeq
They remark that ``it is not known if the invariant ring is generated by the complete weight enumerators of codes''.
This now follows immediately from our main theorem.

Incidentally, the nonzero coefficients of the Molien series in \eqn{EqB1}
form sequence A051462 in \cite{OEIS}, where the reader will find references to both \cite{BDHO99} and \cite{cliff2}.
A great many Molien series arise in studying self-dual codes\footnote{The index to \cite{cliff2} lists the sequence numbers for over 100 
such Molien series.}, and
\cite{OEIS} provides a convenient way to keep track of them.

\section{Example: Hermitian self-dual codes over $\FF_9$}

The form ring for Hermitian self-dual codes over $\FF_9$
 is $\rho = (R=V=\FF_9 , \beta , \Phi )$,
where
\begin{eqnarray*}
\beta (v,w) & = & \frac{1}{3}
\Tr_{ \FF_9 / \FF_3 } ( v \bar{w} ) \\ [+.15in]
\Phi & = & \{ \beta (av, v) , ~a \in \FF_9 \} \\ [+.15in]
\FF_9 & = & \{0,1, \alpha , \alpha^2 , \ldots , \alpha^7 \} \,,
\end{eqnarray*}
with $\alpha^2 + \alpha =1$, $\alpha^4 = -1$.

\vspace*{+.1in}
\noindent
Generators for $\sC (\rho )$ are:
$$
\rho ( \alpha ) = \left[
\begin{array}{ccccccccc}
1 & 0 & 0 & 0 & 0 & 0 & 0 & 0 & 0 \\
0 & 0 & 0 & 0 & 0 & 0 & 0 & 0 & 1 \\
0 & 1 & 0 & 0 & 0 & 0 & 0 & 0 & 0 \\
0 & 0 & 1 & 0 & 0 & 0 & 0 & 0 & 0 \\
0 & 0 & 0 & 1 & 0 & 0 & 0 & 0 & 0 \\
0 & 0 & 0 & 0 & 1 & 0 & 0 & 0 & 0 \\
0 & 0 & 0 & 0 & 0 & 1 & 0 & 0 & 0 \\
0 & 0 & 0 & 0 & 0 & 0 & 1 & 0 & 0 \\
0 & 0 & 0 & 0 & 0 & 0 & 0 & 1 & 0
\end{array}
\right] =: M_1 \,.
$$
The isotropic conditions are:
$$\Phi = \{ \phi (a) : a \in \FF_9 \} \,,$$
\begin{eqnarray*}
\phi (a) (v) & = & \frac{1}{3} \,\Tr ( av \bar{v} ) = \frac{1}{3} \, \Tr (av^4 ) \\
& = & \frac{1}{3} (av^4 + a^3 v^4 ) = \frac{1}{3} ( a+a^3 ) v^4 \,.
\end{eqnarray*}
Take $a= \alpha$, $\alpha + \alpha^3 = -1$. Then
$$\rho ( \phi (\alpha)):~
x_v \mapsto e^{2 \pi i \cdot \frac{-v^4}{3}} x_v
$$
giving the matrix
$$M_2 := {\rm diag} \{1, \omega^2 , \omega , \omega^2 , \omega, \omega^2 , \omega , \omega^2 , \omega \} \,,
$$
where $\omega = e^{ 2 \pi i/3}$.

The MacWilliams transformation:
$$\mbox{idempotent} ~ \iota =1$$
$$x_v \mapsto \frac{1}{\sqrt{9}} \sum_{w \in \FF_9} e^{2 \pi i \frac{1}{3} \, \Tr (\alpha v \bar{w} )} x_w
$$
giving the matrix
$$M_3 := \left[
\begin{array}{ccccccccc}
1 & 1 & 1 & 1 & 1 & 1 & 1 & 1 & 1 \\
1 & \bar{\omega} & \omega & 1 & \omega & \omega & \bar{\omega} & 1 & \bar{\omega} \\
1 & \omega & \omega & \bar{\omega} & 1 & \bar{\omega} & \bar{\omega} & \omega & 1 \\
1 & 1 & \bar{\omega} & \bar{\omega} & \omega & 1 & \omega & \omega & \bar{\omega} \\
1 & \omega & 1 & \omega & \omega & \bar{\omega} & 1 & \bar{\omega} & \bar{\omega} \\
1 & \omega & \bar{\omega} & 1 & \bar{\omega} & \bar{\omega} & \omega & 1 & \omega \\
1 & \bar{\omega} & \bar{\omega} & \omega & 1 & \omega & \omega & \bar{\omega} & 1 \\
1 & 1 & \omega & \omega & \bar{\omega} & 1 & \bar{\omega} & \bar{\omega} & \omega \\
1 & \bar{\omega} & 1 & \bar{\omega} & \bar{\omega} & \omega & 1 & \omega & \omega \\
\end{array}
\right] \,.
$$
Then the Clifford-Weil group is
$$\sC( \rho ) = \langle M_1 , M_2, M_3 \rangle \,,$$
a nine-dimensional group of order 192.

The Molien series for this group is
$$\frac{1+3t^4 + 24t^6 + 74t^8 + 156t^{10} + \cdots + 989 t^{20} + \cdots + t^{38}}{(1-t^2 )^2 (1-t^4 )^2 (1-t^6 )^3 (1-t^8 ) (1-t^{12} )}
$$
\paragraph{Remarks}
\begin{itemize}
\item
The coefficients of the Taylor series expansion form sequence A092354 in \cite{OEIS}.
\item
There are at least 6912 secondary invariants (set $t=1$ in numerator).
\item
This complexity is typical of most groups --- see Huffman and Sloane \cite{HuSl}.
\item
This ring is spanned by cwe's of codes, by our main theorem.
\item
It would be hopeless to try to find a corresponding set of codes!
\item
To get the Hamming weight enumerator theorem, we cannot simply identify
$x_1 = x_\alpha = \cdots = x_{\alpha^7} = x$ (this
fails because $M_2$ does not act nicely, and if we ignore the generator $M_2$ the resulting ring has Molien series
$1/ (1-t^2 )^2$, which is wrong) --- this is what we call an ``illegal symmetrization''.
\item
The correct way to obtain the Hamming weight enumerator theorem is first to divide up the elements of $\FF_9$ into three orbits
$\{0\}$, $\{1, \alpha^2 , \alpha^4 , \alpha^6 \}$ (which square to 1) and
$\{ \alpha , \alpha^3 , \alpha^5 , \alpha^7 \}$ (which square to $-1$).
The generators now collapse nicely, to
$$
\tilde{M}_1 = \left[ \begin{array}{ccc}
1&0&0 \\
0&0&1 \\
0&1&0
\end{array}
\right] \, ,~~\tilde{M}_2 = \left[\begin{array}{ccc}
1&0&0 \\
0 & \omega & 0 \\
0 & 0 & \omega^2
\end{array}
\right]\,,~~
\tilde{M}_3 = \frac{1}{3} \left[
\begin{array}{ccc}
1&1&1 \\
4&1&-2 \\
4&-2&1
\end{array}\right] \,,
$$
generating a group of order 48 with Molien series
$$\frac{1}{(1-t^2 ) (1-t^4 ) (1-t^5)} \,.$$
Codes that correspond to the terms in the denominator can
be taken to be:
$$[1 ~ \alpha ] \,,~~
\left[ \begin{array}{cccc}
1&1&1&0 \\
0&1&2&1
\end{array}\right] \,, ~~
\left[
\begin{array}{cccccc}
1&1&1&1&1&1 \\
1&1&1&0&0&0 \\
0& \alpha & 2 \alpha & 0 & 1 & 2
\end{array}
\right] \,.
$$
If their Hamming weight enumerators are denoted by $f_2$, $f_4$, $f_6$ respectively, then the ring of Hamming weight enumerators is
$$\CC [f_2, f_4] \oplus f_6 \CC [f_2, f_4] \,.$$
This is not the ring of invariants of any finite group of $2 \times 2$ matrices.
\end{itemize}

\section{Higher-genus weight enumerators}

To handle
higher-genus or multiple weight enumerators we use tensor products, as mentioned in \S12,
and Morita theory.
The form ring for genus-$m$ weight enumerators is
$$\rho \otimes R^m = \Mat_m (\rho ) := (\Mat_m (R) , ~V\otimes R^m, ~\beta^{(m)} , ~
\Phi_m ) \,.
$$

\paragraph{Theorem.}
{\em $(1)$ The space of homogeneous invariants of degree $N$ of the corresponding Clifford-Weil group $\sC_m (\rho )$ is spanned by the genus-$m$ weight enumerators $\mbox{cwe}_m (C)$, where $C$ ranges over a set of permutation
representatives of codes of Type $\rho$ and length $N$.
$(2)$~If every length $N$ code of Type $\rho$ is generated by at most $m$
elements, then these genus-$m$ weight enumerators are a basis for the space of homogeneous invariants of degree $N$.
}

\paragraph{Corollary.}
{\em The Molien series of $\sC_m (\rho )$, $Mol_{\sC_m (\rho )} (t)$, converges monotonically as $m$ increases:
$$\lim_{m \to \infty} ~Mol_{\sC_m (\rho )} (t) = \sum_{N=0}^\infty \nu_N t^N \,,$$
where $\nu_N$ is the number of permutation-equivalence classes of codes
of Type $\rho$ and length $N$.
}

\paragraph{Example:}
Binary self-dual (or Type $2_{\rm I}$) codes.
The order of $\sC_m (\rho )$ and the Molien series for genera 1 to 4 are as follows:

\noindent
Genus 1: $|\sC_1 | = 16$ (Gleason \cite{Gle70}):
$$\frac{1}{(1-t^2) (1-t^8 )}$$

\noindent
Genus 2: $|\sC_2 | = 2304$ (see \cite{MMS72}):
$$\frac{1+t^{18}}{(1-t^2 ) (1-t^8 ) (1-t^{12} ) (1-t^{24} )}$$

\noindent
Genus 3: $| \sC_3 | = 5160960$ (see \cite{cliff1}):
$$\frac{\mbox{degree} ~154}{(1-t^2 ) (1-t^{12} ) \cdots (1-t^{40} )}$$

\noindent --- there are at least 720 secondary invariants

\noindent
Genus 4: $| \sC_4 | = 1783 627 77600$ (see Oura \cite{Our97})
$$\frac{\mbox{degree ~504}}{(1-t^2 ) \cdots (1-t^{120})}$$

\noindent --- there are over $10^{10}$ secondary invariants

The convergence of the Molien series mentioned in the above Corollary can be seen in the following table, which gives the initial terms of the expansion
of the Molien series for genera 1--5:
$$
\begin{array}{c|cccccccccc}
m \setminus N & 0 & 2 & 4 & 6 & 8 & 10 & 12 & 14 & 16 & \cdots \\ \hline
1 & 1 & 1 & 1 & 1 & 2 & 2 & 2 & 2 & 3 \\
2 & 1 & 1 & 1 & 1 & 2 & 2 & 3^\ast & 3 & 4 \\
3 & 1 & 1 & 1 & 1 & 2 & 2 & 3 & 4 & 6 \\
4 & 1 & 1 & 1 & 1 & 2 & 2 & 3 & 4 & 7 \\
5 & 1 & 1 & 1 & 1 & 2 & 2 & 3 & 4 & 7
\end{array}
$$
$$
\begin{array}{l}
^\ast\mbox{\footnotesize {This entry 3 corresponds to the fact that the biweight enumerators of the three codes $i_2^6$, $h_8 i_2^2$ and $d_{12}^+$}} \\ [-.05in]
~\mbox{\footnotesize {are linearly independent.}}
\end{array}
$$

\vspace*{+.15in}
Incidentally, the 8-dimensional group $\sC_3 (\rho )$ of order 5160960 is the group whose magical emergence from the computer --- leading to the astonishing coincidence mentioned in \S8 --- indirectly led to our writing the book.

\section{There is no time to mention: }

\begin{itemize}
\item
Our new construction for the Barnes-Wall lattices
as lattices over $\ZZ [\sqrt{2}]$ whose automorphism groups are
the Clifford-Weil groups $\sC_m (2_{\rm I} )$ (Chapter 6).
\item
The theorem that the automorphism group of the genus-$m$ weight enumerator of any Type $2_{\rm I}$ code that is not generated by codewords of weight 2 is
the Clifford-Weil group $\sC_m (2_{\rm I} )$.
There is an analogous assertion for doubly-even or Type $2_{{\rm II}}$ codes.
(Chapter 6)
\item
The generalizations to maximal self-orthogonal codes (Chapter 10).
\item
Quantum codes (Chapter 13).
\item
The extensive tables giving the classification of all codes and of extremal
codes of modest lengths (Chapters 11, 12).
\item
Applications to spherical designs (Chapters 5,6).
\item
``Closed codes'': What definition of duality guarantees that $C^{\perp \perp} = C$?
\item
Our attempts at generalizing the theory to handle self-dual {\em lattices}.
\end{itemize}

\section{Finally, the new list of Types}

Chapter 2 of the book ends with a list of the principal Types and the sections in which they are discussed.
To entice the reader, but without giving any further details,
here is that list:
$$
\begin{array}{|c|c|} \hline
2_{\rm I} {\mbox{~(The~old~Type~I)}} & 4^{\ZZ} {\mbox{~(Codes over~}}\ZZ/4\ZZ) \\ [+.1in]
2_{{\rm II}} {\mbox{~(The~old~Type~II)}} & 4_{{\rm II}}^{\ZZ} \\ [+.1in]
2_S & m^{\ZZ} \\ [+.1in]
2^{\mbox{lin}}, 2_1^{\mbox{lin}} & m_1^{\ZZ} \\ [+.1in]
2_{1'}^{\mbox{lin}}, 2_{1,1'}^{\mbox{lin}} & m_{{\rm II}}^{\ZZ} \\ [+.1in]
4^E {\mbox{~(The~old~Type~IV)}} & m_{{\rm II},1}^{\ZZ} \\ [+.1in]
4_{{\rm II}}^E & m_S^{\ZZ} \\ [+.1in]
q^E \mbox{(even)} & GR(p^e, f)^E \\ [+.1in]
q_{{\rm II}}^E & GR(p^e, f)_{1}^E \\ [+.1in]
3 {\mbox{~(The~old~Type~III)}} & GR(p^e, f)_{p^s}^E \\ [+.1in]
q^E \mbox{(odd)} & GR(2^e, f)_{2^s}^E \\ [+.1in]
q_1^E \mbox{(odd)} & GR (2^e, f)_{{\rm II}}^E \\ [+.1in]
4^H {\mbox{~(The~old~Type~IV)}} & GF (2^e, f)_{{\rm II},2^s}^E \\ [+.1in]
q^H & GR(p^e, f)^H \\ [+.1in]
q_1^H & GR(p^e , f)_{p^s}^H \\ [+.1in]
4^{H+} & GR (p^e , f)^{H+} \\ [+.1in]
4_{{\rm II}}^{H+} & GR (p^e, f)_{p^s}^{H+} \\ [+.1in]
q^{H +} \mbox{(even)} & \ZZ_p  {\mbox{~(Codes over~the~$p$-adic~integers)}} \\ [+.1in]
q_1^{H+} \mbox{(even)} & \FF_{q^2} + \FF_{q^2} u \\ [+.1in]
q_{{\rm II}}^{H+} \mbox{(even)} & \\ [+.1in]
q_{{\rm II},1}^{H+} \mbox{(even)} & \\ [+.1in]
q^{H+} \mbox{(odd)} & \\ [+.1in]
q_1^{H+} \mbox{(odd)} & \\ [+.1in]
q^{\mbox{lin}} , q_1^{\mbox{lin}} & \\ [+.1in]
q_{1'}^{\mbox{lin}}, q_{1,1'}^{\mbox{lin}} & \\ \hline
\end{array}
$$


\end{document}